# Mittag - Leffler function distribution - A new generalization of hyper-Poisson distribution

S. Chakraborty[a] and S. H. Ong[b]

[a]*Department of Statistics, Dibrugarh University, Dibrugarh - 786004, India*

[b]*Institute of Mathematical Sciences, University of Malaya, 50603 Kuala Lumpur, Malaysia*

**Abstract**
In this paper a new generalization of the hyper-Poisson distribution is proposed using the Mittag-Leffler function. The hyper-Poisson, displaced Poisson, Poisson and geometric distributions among others are seen as particular cases. This Mittag-Leffler function distribution (MLFD) belongs to the generalized hypergeometric and generalized power series families and also arises as weighted Poison distributions. MLFD is a flexible distribution with varying shapes like non-increasing with unique mode at zero, unimodal with one / two non-zero modes. It can be under, equi or over dispersed. Various distributional properties like recurrence relation for pmf, cumulative distribution function, generating functions, formulae for different type of moments, their recurrence relations, index of dispersion, its classification, log-concavity, reliability properties like survival, increasing failure rate, unimodality, and stochastic ordering with respect to hyper-Poisson distribution have been discussed. The distribution has been found to fare well when compared with the hyper-Poisson distributions in its suitability in empirical modeling of differently dispersed count data. It is therefore expected that proposed MLFD with its interesting features and flexibility, will be a useful addition as a model for count data.

## 1 Introduction

The Poisson distribution is a popular model for count data. However its use is restricted by the equality of its mean and variance (equi-dispersion). Many models with the ability to represent under, equi and over dispersion have been proposed in the research literature to overcome this restriction. Notable among these distributions are the hyper-Poisson (HP) of Bardwell and Crow [5], generalized Poisson of Consul [8], double-Poisson of Efron et al.[12], Poisson polynomial of Cameron and Johansson [6], weighted Poisson of Castillo and Casany [7] and COM-Poisson of Conway and Maxwell [9] (see also Shmueli et al., [59]).

Of these models the HP distribution which can handle both over- and under-dispersion was first proposed by Bardwell and Crow [5] and Crow and Bardwell [10]. The probability mass function (pmf) of the HP distribution is given by

$$P(X = k) = \frac{\Gamma(\beta)}{\Gamma(k+\beta)} \frac{\lambda^k}{\phi(1,\beta;\lambda)}, \; k = 0,1,2,\cdots; \lambda > 0 \qquad (1)$$

where $\phi(1,\beta;\lambda) = \sum_{k=0}^{\infty} \frac{(1)_k}{(\beta)_k} \frac{\lambda^k}{k!}, (\beta)_k = \beta(\beta+1)\cdots(\beta+k-1)$ is the confluent hypergeometric function. The pmf has a simple recurrence relation given by

$$(\beta+k)P(X=k+1) = \lambda P(X=k), \quad k=1,2,3,\ldots.$$

The probability generating function (pgf) is

$$P(s) = {}_1\phi_1(1,\beta;\lambda s)/{}_1\phi_1(1,\beta;s).$$

Staff [62] studied displaced Poisson distribution which is HP distribution with the parameter $\beta$ restricted to be a positive integer. The case when $\beta$ is negative was investigated later by Staff [63]. The HP distribution drew attention of many researchers of late. Kemp [29] dealt with a q-analogue of the distribution and Ahmad [1] proposed a Conway-Maxwell-HP distribution. Roohi and Ahmad [53, 54] investigated moments of the HP distribution. Kumar and Nair [34, 35, 36, and 37] studied various extensions and alternatives of the HP distribution. Castillo et al. [7]) studied a HP Regression Model for over-dispersed and under-dispersed count data. Best [3] and Antic et al [2] considered the HP distribution in word length and text length research. Khazraee et al. [31] investigated the application of HP generalized linear model for analyzing motor vehicle crashes.

In this paper we wish to propose a new generalization of the HP distribution by replacing $\Gamma(k+\beta)$ in (1) with $\Gamma(\alpha k + \beta), \alpha > 0$ and the normalization constant becomes $E_{\alpha,\beta}(\lambda)$ which is the generalized Mittag-Leffler function defined by

$$E_{\alpha,\beta}(z) = \sum_{k=0}^{\infty} z^k / \{\Gamma(\alpha k + \beta)\}. \tag{2}$$

Consequently the proposed distribution is called the Mittag-Leffler function distribution (MLFD). The extra parameter $\alpha$ adds flexibility to the model but retains computational tractability since computation of $E_{\alpha,\beta}(\lambda)$ does not pose a problem due to many software packages (example MATLAB) which offered routines for its quick computation The MLFD is shown to be log-concave and this confers a number of attractive properties for modeling and inference; see Walther [66] for a good review of statistical modeling and inference with log-concave distributions. The proposed MLFD should not be confused with a class of discrete Mittag-Leffler distribution proposed by Pillai and Jayakumar [51]. It is pertinent to give a brief review of some developments in statistical models involving the Mittag-Leffler function.

The Mittag-Leffler function (with $\beta = 1$ in (2) and is denoted by $E_\alpha(z)$) was first introduced by Swedish mathematician Gosta Mittag-Leffler [43, 44]) and it arises as the solution of a fractional differential equation. This function and its many extended versions were studied by many mathematicians over the years. Haubold, Mathai and Saxena [21] have given a good survey on the Mittag-Leffler function. Recently, this function has also been explored for applications in statistics. Pillai [50] shown that $1 - E_\alpha(-z^\alpha), 0 < \alpha \leq 1$ are valid cumulative distribution functions (cdf) and named it as Mittag-Leffler distribution with pdf and cdf respectively given by

$$F(x;\alpha) = \sum_{k=1}^{\infty} (-1)^{k-1} x^{k\alpha} / \{\Gamma(\alpha k + 1), x > 0, 0 < \alpha \leq 1 \text{ and}$$

$$f(x;\alpha) = \sum_{k=1}^{\infty}(-1)^{k-1}(k\alpha)x^{k\alpha-1}/\{\Gamma(\alpha k+1), x>0, 0<\alpha \leq 1 \qquad (3)$$

Since for $\alpha = 1$ this distribution reduces to exponential distribution with mean 1, it can be treated as a generalization of the exponential distribution. He studied different properties of this distribution. Jayakumar and Pillai [23], Jose and Pillai [26], Lin [39], Jayakumar [24], Joes et al. [27] studied different aspects of this distribution.

Pillai and Jayakumar [51] proposed a class of discrete Mittag-Leffler distribution (DML) having probability generating function (pgf) $P(z) = E(X^z) = 1/[1+c(1-z)^{\alpha}]$. The DML distribution arises as a mixture of Poisson with parameter $\theta \lambda$ with $\theta$ is a constant and $\lambda$ following the Mittag-Leffler distribution in (2). They have studied different properties of DML and gave a probabilistic derivation and application in a first order auto regressive discrete process. DML is also a particular case of the discrete Linnik distribution [11].

Jose and Abraham [28] introduced another discrete distribution based on Mittag-Leffler function. Their distribution is a kind of generalization of Poisson distribution in the sense that is arises when the exponential waiting time distribution in the usual Poisson process is replaced by the Mittag-Leffler distribution. The pmf of their distribution is

$$P(X=k) = \sum_{i=k}^{\infty}\binom{i}{k}(-1)^{(i-k)}z^{i\alpha}/\{\Gamma(\alpha i+1), k=0,1,\cdots; 0<\alpha \leq 1 \qquad (3)$$

In this article we have taken a completely different route to propose a discrete distribution based on the Mittag-Leffler function. The proposed distribution which also arises from a queuing theory set up is simple and extremely flexible in its shape and modality, it can model under-dispersed, equi-dispersed and over-dispersed count data. Section 2 defines the Mittag-Leffler function distribution (MLFD) and basic structural properties are given. MLFD as a distribution in a queuing system is given in section 3. Reliability and stochastic ordering properties are discussed in section 4. Section 5 considers the computation of the generalized Mittag-Leffler function while section 6 gives examples of applications of MLFD. The conclusion is given in section 7.

## 2 Mittag-Leffler Function Distribution: Definition and Properties

**Definition 1.** A discrete random variable $X$ is said to follow the Mittag-Leffler function distribution (MLFD) if its pmf is given by

$$P(X=k) = \lambda^k/\{\Gamma(\alpha k+\beta)E_{\alpha,\beta}(\lambda)\}, k=0,1,2,\cdots; \lambda,\alpha,\beta>0 \qquad (4)$$

where $E_{\alpha,\beta}(z) = \sum_{k=0}^{\infty} z^k/\{\Gamma(\alpha k+\beta)$ is the generalized Mittag-Leffler function. The distribution henceforth will be denoted by MLFD $(\lambda,\alpha,\beta)$.

### 2.1 Probability recurrence relation
The MLFD $(\lambda,\alpha,\beta)$ pmf in (4) has a simple recurrence relation given by

$$\Gamma(\alpha k+\alpha+\beta)P(X=k+1) = \lambda\, \Gamma(\alpha k+\beta)P(X=k), k\geq 1 \qquad (5)$$

with $P(X=0) = 1/\{\Gamma(\beta)E_{\alpha,\beta}(\lambda)\}$.

When $\alpha$ a positive integer, (5) can be expressed as $(\alpha k+\beta)_{\alpha} P(X=k+1) = \lambda P(X=k)$,

where $(y)_\alpha = y(y+1)...(y+\alpha-1)$.

The distribution exhibits longish tail for $0 < \alpha < 1$ as the ratio of successive probabilities varies slowly (this corresponds to over dispersion) as *k* tends to infinity while for $\alpha \geq 1$ this ratio tends to zero faster implying presence of a Poisson type tail.

The above recurrence relation facilitates easy computation of the probabilities. The computation of the normalizing constant $E_{\alpha,\beta}(\lambda)$ is only required for $P(X = 0)$.

Note that the recurrence relation (or the difference equation) in (5) reduces respectively to that of Hyper-Poisson [5] and displaced Poisson distribution [62] for $\alpha = 1$ and $\alpha = 1, \beta \geq 0$ an integer (further discussed in see section 2.4).

## *2.2 Shapes of probability mass function*

The pmf of $\text{MLFD}(\lambda, \alpha, \beta)$ is plotted for a number of combinations of parameters to study the different shapes of the distribution.

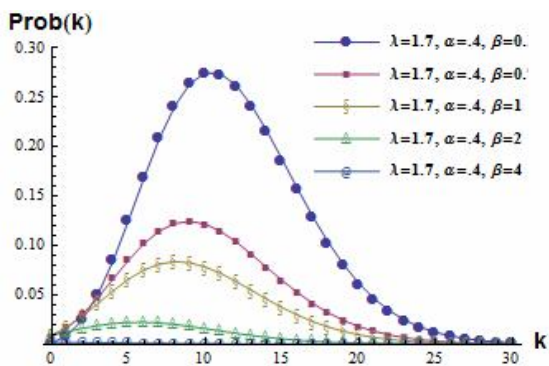

Fig. 1

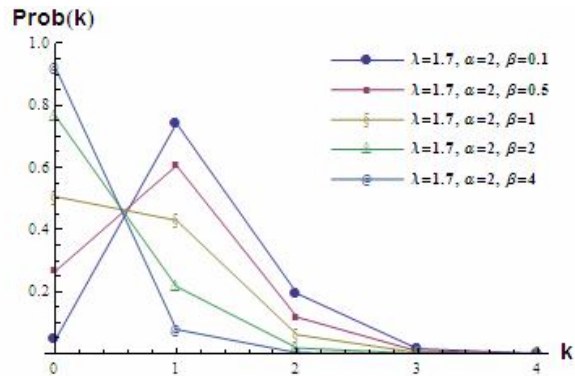

Fig. 2

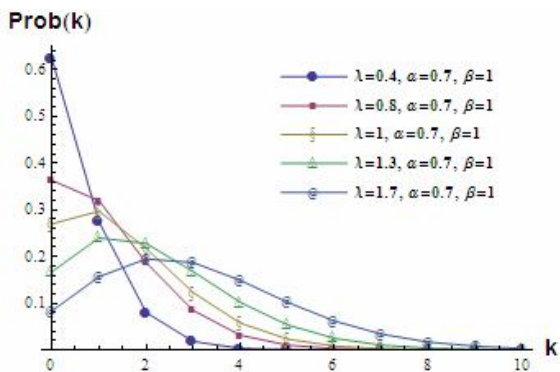

Fig. 3

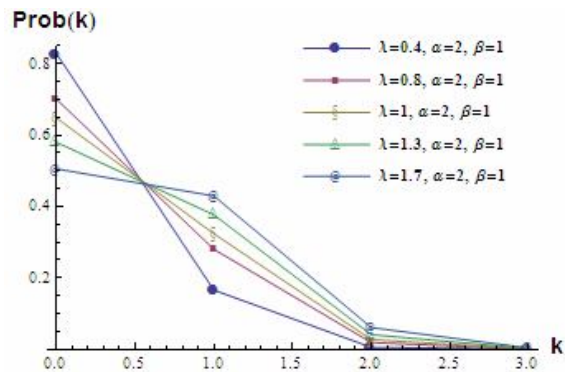

Fig. 4

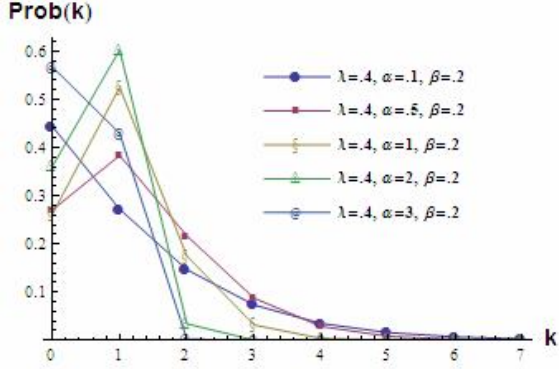
Fig. 5

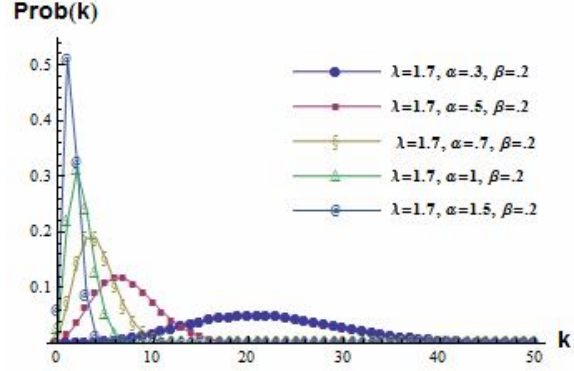
Fig. 6

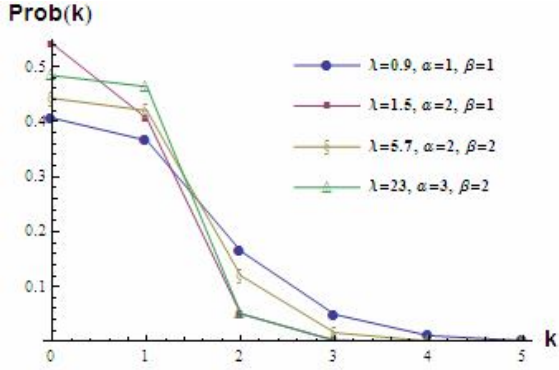
Fig. 7

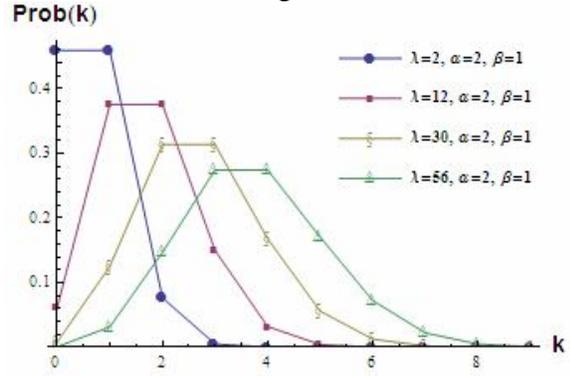
Fig. 8

From the plots of the pmf it is seen that the distribution can be unimodal with nonzero mode (see Fig. 1) or it can have two nonzero modes (see Fig. 8) or non-increasing with the mode at 0 (see Fig. 7) (see section 4 for further discussion on modes).

### *Dependence of the pmf on the parameters*:

The proportion of zeros given by $P(X=0) = 1/\{\Gamma(\beta)E_{\alpha,\beta}(\lambda)\}$, increases with increase in $\beta$ (see Fig. 2) when $(\lambda,\alpha)$ are fixed, increases with decrease in $\lambda$ (see Fig. 3, 4) for fixed $(\alpha,\beta)$ and increases with increase in $\alpha$ (see Fig. 5 and 6) for fixed $(\lambda,\beta)$.

Whereas when $\beta \to 0^+$ the proportion zero decreases (can be seen in the Fig. 2) with the pmf proportional to $\lambda^k /\{\Gamma(\alpha k)$, $k=0,1,2,\cdots$.

### *2.3 Cumulative distribution function and generating functions*

The cumulative distribution function (cdf) of MLFD$(\lambda,\alpha,\beta)$ is seen to be

$$P(X \leq r) = 1 - \lambda^{r+1}[E_{\alpha,\beta+(r+1)\alpha}(\lambda)/E_{\alpha,\beta}(\lambda)]$$

by using the known relation $\lambda^r E_{\alpha,\beta+r\alpha}(\lambda) = E_{\alpha,\beta}(\lambda) - \sum_{k=0}^{r-1} z^k /\{\Gamma(\alpha k+\beta)\}$ (Haubold et al. [21]).

The pgf of MLFD$(\lambda,\alpha,\beta)$ is given in terms of $E_{\alpha,\beta}(\lambda)$ as:

$$P(s) = E(s^X) = E_{\alpha,\beta}(\lambda s)/E_{\alpha,\beta}(\lambda) \quad , \quad 0 < \lambda s \tag{6}$$

The moment generating function (mgf) $M_X(s)$ and the factorial moment generating function (fmgf) are respectively given by

$$P(e^t) = E_{\alpha,\beta}(\lambda e^s)/E_{\alpha,\beta}(\lambda) \text{ and } P(1+t) = E_{\alpha,\beta}(\lambda(1+t))/E_{\alpha,\beta}(\lambda)$$

## 2.4 Related distributions and connections with other families of distributions
### 2.4.1 Particular cases of MLFD $(\lambda, \alpha, \beta)$

The MLFD $(\lambda, \alpha, \beta)$ includes a number of well-known distributions as particular case:

(i) When $\alpha = \beta = 1$, MLFD $(\lambda, \alpha, \beta)$ reduces to the **Poisson** distribution with parameter $\lambda$.

(ii) When $\alpha = 0$, $\beta (\geq 0)$ MLFD $(\lambda, \alpha, \beta)$ becomes **geometric** distribution with parameter $\lambda$ provided $|\lambda| < 1$.

Since for $\alpha \to 0^+$,

$$\lim_{\alpha \to 0^+} E_{\alpha,\beta}(\lambda) = \sum_{k=0}^{\infty} \frac{\lambda^k}{\Gamma(\alpha k + \beta)} \to \sum_{k=0}^{\infty} \frac{\lambda^k}{\Gamma\beta} = \frac{1}{\Gamma\beta(1-\lambda)}, |\lambda| < 1. \text{(Hanneken et al. [20])}$$

(iii) When $\alpha = 1$, $\beta (\geq 0)$ MLFD $(\lambda, \alpha, \beta)$ reduces to the HP $(\lambda, \beta)$ distribution (Bardwell and Crow [5], p. 200 Johnson et al.[25])

*Proof*: $P(X = k) = \frac{\lambda^k}{\Gamma(k+\beta)E_{1,\beta}(\lambda)}$, $k = 0, 1, 2, \cdots; \lambda > 0$ where

$$E_{1,\beta}(\lambda) = \sum_{k=0}^{\infty} \lambda^k / \Gamma(k+\beta) = \sum_{k=0}^{\infty} \lambda^k / \{(\beta)_k \Gamma(\beta)\} = \frac{1}{\Gamma(\beta)} \sum_{k=0}^{\infty} \frac{(1)_k}{(\beta)_k} \frac{\lambda^k}{k!} = \frac{\phi(1,\beta;\lambda)}{\Gamma(\beta)},$$

where $\phi(1, \beta; \lambda)$ is the confluent hypergeometric function. Hence

$$P(X = k) = \frac{\Gamma(\beta)}{\Gamma(k+\beta)} \frac{\lambda^k}{\phi(1,\beta;\lambda)}, k = 0, 1, 2, \cdots; \lambda > 0.$$

An alternative form of the above pmf can be seen as

$$P(X = k) = \frac{\Gamma(\beta-1)}{\Gamma(k+\beta)} \frac{e^{-\lambda} \lambda^{k+\beta-1}}{\gamma(\beta-1,\lambda)}, k = 0, 1, 2, \cdots; \lambda > 0$$

since $E_{1,\beta}(\lambda) = \frac{\lambda^{1-\beta} e^\lambda \gamma(\beta-1,\lambda)}{\Gamma(\beta-1)}$ (see Simon, [60]), where $\gamma(u,\lambda) = \int_0^\lambda e^{-s} s^{u-1} ds$ is incomplete gamma function.

(iv) When $\alpha = 1$ and $\beta (= t+1)$ is a positive integer, MLFD $(\lambda, \alpha, \beta)$ reduces to the displaced Poisson distribution (see Staff [62], p. 200 of Johnson et al., [25]) with parameter $\lambda$ and $t$.

In addition the following new distributions involving hyperbolic, and error functions are also seen as particular cases.

(v) When $\alpha = 2$ and $\beta = 2$, MLFD $(\lambda, \alpha, \beta)$ reduces to a new discrete distribution with parameter $\lambda$ and pmf

$$P(X=k) = \frac{\lambda^{k+(1/2)}}{\Gamma 2(k+1)} \frac{1}{\sinh(\sqrt{\lambda})} = \frac{(\sqrt{\lambda})^{2k+1}}{(2k+1)!} \frac{1}{\sinh(\sqrt{\lambda})}, \quad k=0,1,2,\cdots; \lambda > 0$$

since $E_{2,2}(\lambda) = \sinh(\sqrt{\lambda})/\sqrt{\lambda}$. (Haubold et al., [21]).

(vi) When $\alpha = 2$ and $\beta = 1$, MLFD$(\lambda, \alpha, \beta)$ reduces to a new distribution with parameter $\lambda$ and pmf

$$P(X=k) = \frac{\lambda^k}{\Gamma(2k+1)} \frac{1}{\cosh(\sqrt{\lambda})} = \frac{(\sqrt{\lambda})^{2k}}{(2k)!} \frac{1}{\cosh(\sqrt{\lambda})}, \quad k=0,1,2,\cdots; \lambda > 0$$

since $E_{2,1}(\lambda) = \cosh(\sqrt{\lambda})$. (Haubold et al., [21]).

(vii) When $\alpha = 1/2$ and $\beta = 1$, MLFD$(\lambda, \alpha, \beta)$ reduces to a new distribution with parameter $\lambda$ and pmf

$$P(X=k) = \frac{\exp(-\lambda^2)\lambda^k}{(k/2)!\ erfc(-\lambda)}, \quad k=0,1,2,\cdots; \lambda > 0$$

since $E_{1/2,1}(\sqrt{\lambda}) = \exp[\lambda]erfc(-\sqrt{\lambda})$ (Haubold et al., [21]) where $erfc(\lambda)$ is the complementary error function defined as $erfc(\lambda) = 1 - erf(\lambda) = 1 - \frac{2}{\sqrt{\pi}} \int_0^\lambda \exp(-t^2) dt$.

Also $erfc(\lambda) = 2\left[\frac{1}{\sqrt{2\pi}} \int_0^{-\sqrt{2}\lambda} \exp(-t^2/2) dt\right] = 2\Phi(-\sqrt{2}\lambda)$, $\Phi(.)$ being the cdf of standard normal distribution.

**Remark 1.**
a. For $0 \leq \alpha \leq 1$, the MLFD$(\lambda, \alpha, \beta)$ can be viewed as a *continuous bridge* between geometric ($\alpha = 0$) and hyper-Poisson ($\alpha = 1$) distributions in the range of the parameter $\alpha$.
b. For $0 \leq \alpha \leq 1$, the MLFD$(\lambda, \alpha, 1)$ can be viewed as a *continuous bridge* between geometric ($\alpha = 0$) and Poisson ($\alpha = 1$) distributions in the range of the parameter $\alpha$. This property is also shared by the Com-Poisson distribution [9].

### *2.4.2 MLFD as weighted distributions*
i. MLFD$(\lambda, \alpha, \beta)$ as weighted MLFD$(\lambda, \alpha, 1)$: If $X \sim$ MLFD$(\lambda, \alpha, 1)$ having pmf

$$P(X=k) = \frac{\lambda^k}{\Gamma(\alpha k+1) E_\alpha(\lambda)}, \quad k=0,1,2,\cdots; \lambda,>0,$$

then for integer $\beta$ it can be shown that weighted distribution with weight function $1/(\alpha k+1)_\beta$, gives the pmf of MLFD$(\lambda, \alpha, \beta)$.

ii. MLFD$(\lambda, \alpha, \beta)$ as weighted HP$(\lambda, \beta)$ (i.e. MLFD$(\lambda, 1, \beta)$):
If $X \sim$ MLFD$(\lambda, 1, \beta)$ having pmf

$$P(X=k) = \frac{\Gamma(\beta)}{\Gamma(k+\beta)} \frac{\lambda^k}{\phi(1,\beta;\lambda)}, \quad k=0,1,2,\cdots; \lambda,>0,$$

Then for integer $\alpha > 1$, it can be shown that weighted distribution with weight function $1/(k+\beta+1)_{(\alpha-1)k}$, gives the pmf of MLFD $(\lambda, \alpha, \beta)$.

iii. MLFD $(\lambda, \alpha, \beta)$ as weighted Poisson: If $X \sim$ Poisson $(\lambda)$ having pmf
$$P(X = k) = e^{-\lambda}\lambda^k / k!, \; k = 0,1,2,\cdots; \lambda > 0,$$
Then for integer $\alpha$ and $\beta$ it can be shown that weighted distribution with weight function $1/(k+1)_{(\alpha-1)k+\beta-1}$ and $1/(k+1)_{\beta-1}$ respectively gives the pmf of MLFD $(\lambda, \alpha, \beta)$ and MLFD $(\lambda, 1, \beta)$.

**Remark 2**. Since the weight function $1/(k+\beta+1)_{(\alpha-1)k}$ is monotonically decreasing in $k$ for $\alpha > 1$ (see Patil et al., [49]; Ross [55]; Castillo and Casany [7]), MLFD $(\lambda, \alpha, \beta)$ is stochastically smaller than the hyper-Poisson and hence may be seen as a sub-hyper-Poisson distribution when $\alpha > 1$.

### 2.4.3 Limiting distributions of MLFD

i. For $\lambda \to \infty$, MLFD $(\lambda, \alpha, \beta)$ tends to a new distribution with pmf
$$P(X = k) = \frac{\alpha \lambda^{(\beta+k-1)/\alpha} \exp(-\lambda^{1/\alpha})}{\Gamma(\alpha k + \beta)}, \; k = 0,1,2,\cdots; \lambda > 0$$
This follows on using the result that $E_{\alpha,\beta}(\lambda) \to \{\lambda^{(1-\beta)/\alpha}\exp(\lambda^{1/\alpha})\}/\alpha$ for large values of $\lambda$ ( Gerhold [16]). In particular, MLFD $(\lambda, \alpha, 1)$ tends to a new distribution with pmf
$$P(X = k) = \frac{\alpha \lambda^k \exp(-\lambda^{1/\alpha})}{\Gamma(\alpha k + 1)}, \; k = 0,1,2,\cdots; \lambda > 0$$
This follows on using the result that $E_\alpha(\lambda) \to \exp(\lambda^{1/\alpha})/\alpha$ for large values of $\lambda$ ( Erdelyi, [13], p. 208; Gerhold, [16]). Both these approximations are good when $\alpha \in (0, 2]$ (see Simon, [60]).

ii. MLFD $(\lambda, \alpha, \beta)$ degenerates with mass only at zero for either $\alpha \to \infty$ or $\beta \to \infty$ or both and also when $\lambda \to 0^+$.

### 2.4.4 MLFD as member of some families of discrete distributions

i. MLFD $(\lambda, \alpha, \beta)$ is a member of the generalized hypergeometric family (Kemp [30], see Johnson et al. [25] for details). This is clear from equation (6).

ii. MLFD $(\lambda, \alpha, \beta)$ is a member of the generalized power series distribution (Patil [46, 47]) when $\lambda$ is the primary parameter.
The GPSD has the pmf of the form
$P(X = x) = A(x)\lambda^x / f(\lambda)$, $x \in$ a sub set of non-negative integers $A(x) \geq 0$, $\lambda > 0$ so that $f(\lambda)$ is positive, finite and differentiable functions of $\lambda$. In case of MLD $(\lambda, \alpha, \beta)$, $A(x) = 1/\Gamma(\alpha x + \beta)$ and $f(\lambda) = E_{\alpha,\beta}(\lambda)$.

iii. For fixed value of the parameter $\alpha$ and $\beta$ the MLFD $(\lambda, \alpha, \beta)$ is also a member of the exponential family of distribution

To see this consider the likelihood for a set of $n$ independent and identically distributed observations $x_1, x_2, \cdots, x_n$ from MLFD$(\lambda, \alpha, \beta)$ the likelihood function is given by

$$L(x_1, x_2, \cdots, x_n | z, \alpha, \beta) = \lambda^{\sum_{i=1}^{n} x_i} [E_{\alpha,\beta}(\lambda)]^{-n} \prod_{i=1}^{n} \Gamma(\alpha x_i + \beta) \}^{-1}$$

$$= \lambda^{\sum_{i=1}^{n} x_i} e^{-\sum_{i=1}^{n} \log \Gamma(\alpha x_i + \beta)} [E_{\alpha,\beta}(\lambda)]^{-n} = \lambda^{S_1} e^{-S_2} [E_{\alpha,\beta}(\lambda)]^{-n} \quad (7)$$

where $S_1 = \sum_{i=1}^{n} x_i$, $S_2 = \sum_{i=1}^{n} \Gamma(\alpha x_i + \beta)$.

iv. When the parameter $\alpha$ is integer, the pgf of MLFD$(\lambda, \alpha, 1)$ can be represented in terms of confluent hypergeometric function as

$$P(s) = \frac{E_\alpha(\lambda s)}{E_\alpha(\lambda)} = {}_0F_{\alpha-1}\left(; \frac{1}{\alpha}, \frac{2}{\alpha}, \cdots, \frac{\alpha-1}{\alpha}; \frac{\lambda s}{\alpha^\alpha}\right) \bigg/ {}_0F_{\alpha-1}\left(; \frac{1}{\alpha}, \frac{2}{\alpha}, \cdots, \frac{\alpha-1}{\alpha}; \frac{\lambda}{\alpha^\alpha}\right)$$

Thus MLFD$(\lambda, \alpha, 1)$ belongs to generalized hypergeometric family of Kemp [30].

## 2.5 Moments and related results

Denoting $E(X^r) = \mu'_r$, $E(X)_{[r]} = \mu_{[r]}$ and $E[\{X - E(X)\}^r] = \mu_r$, $y_{[\alpha]} = y(y-1)\cdots(y-\alpha+1)$ and using the relation $E(X_{[r]}) = \mu'_r = \frac{d^i}{ds^i} P(s)\big|_{s=1}$ where $P(s)$ is the pgf of MLFD$(\lambda, \alpha, \beta)$ mentioned in section 2.3 and along with a result of derivative of $E_{\alpha,\beta}(\lambda)$ with respect to $\lambda$ given by

$$\frac{d}{d\lambda} E_{\alpha,\beta}(\lambda) = \frac{E_{\alpha,\beta-1}(\lambda) - (\beta-1) E_{\alpha,\beta}(\lambda)}{\alpha \lambda},$$

the following formulae can be derived:

$$\mu'_1 = \lambda \frac{d}{d\lambda} \log[E_{\alpha,\beta}(\lambda)] = \frac{E_{\alpha,\beta-1}(\lambda)}{\alpha E_{\alpha,\beta}(\lambda)} + \frac{1-\beta}{\alpha}, \text{ provided } \alpha > 0 \text{ and } \beta > 1.$$

$$\mu'_2 = \frac{1}{\alpha^2} \frac{E_{\alpha,\beta-2}(\lambda)}{E_{\alpha,\beta}(\lambda)} - \frac{2\beta-3}{\alpha^2} \frac{E_{\alpha,\beta-1}(\lambda)}{E_{\alpha,\beta}(\lambda)} + \left(\frac{\beta-1}{\alpha}\right)^2$$

$$\mu_2 = \frac{1}{\alpha^2} \left\{ \frac{E_{\alpha,\beta-2}(\lambda)}{E_{\alpha,\beta}(\lambda)} - \left(\frac{E_{\alpha,\beta-1}(\lambda)}{E_{\alpha,\beta}(\lambda)}\right)^2 + \frac{E_{\alpha,\beta-1}(\lambda)}{E_{\alpha,\beta}(\lambda)} \right\}, \text{ provided } \alpha > 0 \text{ and } \beta > 2.$$

The variance can also be expressed as

$$\mu_2 = \lambda \frac{d}{d\lambda} \mu'_1 = \lambda \frac{d}{d\lambda} \left[ \lambda \frac{d}{d\lambda} \log[E_{\alpha,\beta}(\lambda)] \right] = \lambda \frac{d}{d\lambda} \log[E_{\alpha,\beta}(\lambda)] + \lambda^2 \frac{d^2}{dp^2} \log[E_{\alpha,\beta}(\lambda)].$$

The above results can alternatively be derived easily by first deriving $E(\alpha X + \beta)_{[r]}, r = 1, 2$ and then simplifying.

***Caution***: In all the above expressions there is restriction on the values of $\beta$. This situation **may** be overcome by using the following relation repeatedly till the conditions are satisfied.

$E_{\alpha,\beta}(\lambda) = \frac{1}{\Gamma(\beta)} + \lambda E_{\alpha,\alpha+\beta}(\lambda)$ (see Erdelyi et al. [13], Hanneken et al., [20]), the gamma function for negative $\beta$ can be computed using the formula (see Fisher and Kilicman, [14])

$$\Gamma(-n) = \begin{cases} -\Gamma(-n+1)/n, \text{ when } n \neq 1,2,\cdots \\ (-1)^n / n!\{\phi(n) - \gamma\}, n = 1,2,\cdots \end{cases}$$

where $\phi(n) = \sum_{i=1}^{n} \frac{1}{i}$ and $\gamma = -\Gamma'(1)$ is the Euler's constant.

*Alternative formulae for moments*:

To avoid the complicacy that might arise in the above formulations alternative simple formulae for moments in terms of generalized ML function of Prabhakar [52] may be presented as follows:

$\mu_1' = \frac{E_{\alpha,\beta}^2(\lambda)}{E_{\alpha,\beta}(\lambda)} - 1$, $\mu_2' = \frac{2E_{\alpha,\beta}^3(\lambda) - 3E_{\alpha,\beta}^2(\lambda)}{E_{\alpha,\beta}(\lambda)} + 1$ etc. where $E_{\alpha,\beta}^\rho(\lambda) = \sum_{k=0}^{\infty} \frac{(\rho)_k}{k!} \frac{\lambda^k}{\Gamma(\alpha k + \beta)}$ is the generalized Mittag Lefller function of Prabhakar [52] where $(y)_\alpha = y(y+1)\ldots(y+\alpha-1)$.

In general, $E[(X+1)_r] = E_{\alpha,\beta}^{r+1}(\lambda)/E_{\alpha,\beta}(\lambda)$.

*General formulae for the $r^{th}$ factorial moment*:

Following equations (6.5) and (6.8) of Janardan [22] the following results can be written easily

$\mu_{[r]} = \lambda^r \frac{E_{\alpha,\beta}^{(i)}(\lambda)}{E_{\alpha,\beta}(\lambda)}$ and $\mu_r' = \lambda^r \frac{E_{\alpha,\beta}(\lambda \mathbf{E})0^r}{E_{\alpha,\beta}(\lambda)} = \sum_{i=0}^{r} \frac{\lambda^i}{E_{\alpha,\beta}(\lambda)} E_{\alpha,\beta}^{(i)}(\lambda) \frac{\Delta^i 0^r}{i!}$,

where $E_{\alpha,\beta}^{(i)}(\lambda)$ is the ith derivative of $E_{\alpha,\beta}(\lambda)$ with respect to $\lambda$ and $\mathbf{E} p(x) \equiv (1+\Delta)p(x) = p(x+1)$ is the *shift operator*, $\Delta p(x) = p(x+1) - p(x)$ is the *forward difference operator* with interval of differencing being equal to '1' and $E^i 0^r = E^i Y^r \big|_{y=0}$.

*Recurrence relations of moments*:

The following recurrence relations hold:

(i) $\mu_{r+1}' = \lambda \frac{d}{dp} \mu_r' + \mu_r' \mu_1'$

(ii) $\mu_{r+1} = \lambda \frac{d}{dp} \mu_r + r\mu_{r-1} \mu_2$

(iii) $\mu_{[r+1]} = \lambda \frac{d}{dp} \mu_{[r]} + \mu_{[r]} \mu_1' - r \mu_{[r]}$

(iv) $E\{\alpha(X-1) + \beta\}_\alpha = \lambda + (\beta - \alpha)_\alpha P(X=0)$

The relations (i) to (iii) can be proved using the general relations for GPSD or by direct manipulation while (iv) follows from difference equation in (5).

*Monotonicity of the mean*:

Since $\mu_2 > 0$ and $\lambda \neq 0$, $\mu_2 = \lambda \frac{d}{d\lambda}\mu_1' > 0$ implies $\frac{d}{d\lambda}\mu_1' > 0$. Hence $\mu_1'$ is a monotonically increasing function of $\lambda$.

*A useful moment identity*:

The following useful moment identity for MLFD$(\lambda, \alpha, \beta)$ can be derived following Kattumannil and Tibilet [32]:

$$E[c(X)(X - \mu_1')] = -\lambda\, E\left[\Delta c(x) \frac{1}{P(x)} \frac{d}{d\lambda} F(x)\right],$$

where $c(x)$ is a real function of $x$, $\Delta c(x) = c(x+1) - c(x)$, $F(x)$ and $p(x)$ are respectively the cdf and pmf of MLFD$(\lambda, \alpha, \beta)$.

Choosing an appropriate form of $c(x)$ one can get different results. For example, it can be easily checked that by taking $c(X) = (X - \mu_1')$ and $c(X) = (X + \mu_1')$ we can recapture the relations

$$\mu_2 = \lambda \frac{d}{d\lambda}\mu_1' \text{ and } \mu_2' = \lambda \frac{d}{d\lambda}\mu_1' + (\mu_1')^2 \text{ respectively.}$$

## 2.6 Index of dispersion

The index of dispersion (ID) is given by ID = Variance / Mean. Contour plots of ID for different choices of parameters $(\alpha, \beta)$ keeping $\lambda = 0.25$ and $\lambda = 5$ fixed are presented in Fig. 9 and 10 to depict the contours of ID. Labels on a given line indicate the value of ID on that line. In these figures, the ID of the region on the left (right) of a given line is more (less) than the ID value on that given line.

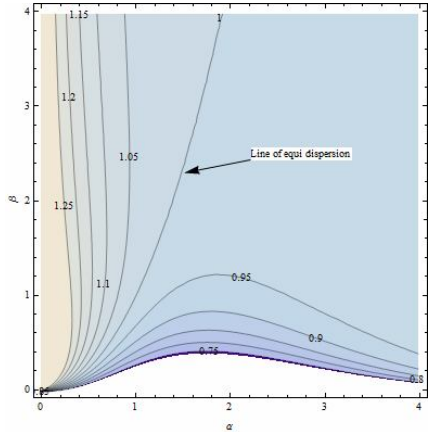
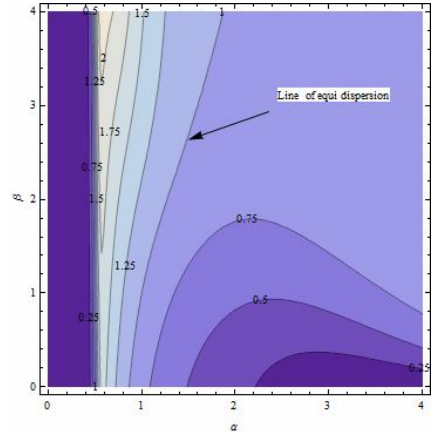

Fig. 9  Fig. 10

From the Fig. 9 and 10, it is obvious that the MLFD$(\lambda, \alpha, \beta)$ is very flexible with respect to the ID and able to cater for under, equi and over dispersion in count data. Interestingly, this family includes a non Poisson distribution with equi-dispersion. We next consider a theorem to state the conditions under which the distribution is equi, under or over dispersed.

**Theorem 1.** The MLFD$(\lambda, \alpha, \beta)$ is

(i) equi-dispersed, if $E_{\alpha,\beta}(\lambda) = \exp[a + b\lambda]$

(ii) over-dispersed, if $E_{\alpha,\beta}(\lambda) = \exp[a + b\lambda + h(\lambda)]$

(iii) under-dispersed, if $E_{\alpha,\beta}(\lambda) = \exp[a + b\lambda + g(\lambda)]$

where $a$ and $b$ are arbitrary constants and $h(\lambda) > 0$ is increasing and $g(\lambda) < 0$ is decreasing function of $\lambda$.

*Proof*: For MLFD $(\lambda, \alpha, \beta)$ the ID is given by

$$\text{ID} = \frac{\mu_2}{\mu_1'} = \frac{1}{\mu_1'}\lambda\frac{d}{d\lambda}\mu_1' = \lambda\frac{d}{d\lambda}\log\mu_1' = \lambda\frac{d}{d\lambda}\log\left[\lambda\frac{d}{d\lambda}\log[E_{\alpha,\beta}(\lambda)]\right]$$

(i) ID = 1. We have

$$\text{ID} = \lambda\left[\frac{d}{d\lambda}\log\lambda + \frac{d}{d\lambda}\log\frac{d}{d\lambda}[\lambda E_{\alpha,\beta}(\lambda)]\right] = 1.$$

That is, $1 + \frac{d}{d\lambda}\log\frac{d}{d\lambda}[\lambda E_{\alpha,\beta}(\lambda)] = 1$ or $\frac{d}{d\lambda}\log\frac{d}{d\lambda}[\lambda E_{\alpha,\beta}(\lambda)] = 0$, since $\lambda \neq 0$.

Hence $\log\frac{d}{d\lambda}[\lambda E_{\alpha,\beta}(\lambda)] = c$.

Thus it follows that $E_{\alpha,\beta}(\lambda) = \exp[c\lambda + b]$, $c\ (>0)$ and $b$ are arbitrary constants.

(ii) Similarly ID > 1 implies

$$E_{\alpha,\beta}(\lambda) = \exp[c\lambda + b + h(\lambda)], \ h(\lambda) > 0 \text{ is increasing function of } \lambda.$$

(iii) ID < 1 implies

$$E_{\alpha,\beta}(\lambda) = \exp[c\lambda + b + g(\lambda)], \ g(\lambda) < 0 \text{ is decreasing function of } \lambda.$$

*Numerical verification of Theorem 1*: We have verified the results of the above theorem by plotting $\log E_{\alpha,\beta}(\lambda)$ against $\lambda$ for some pairs of values of $\alpha$ and $\beta$ taken from line of equi-dispersion in the contour plots in Fig. 9 and 10. Here few representative plots are presented in Fig. 11. It can be easily checked that $Log[E_{\alpha,\beta}(\lambda)]$ is a linear function of $\lambda$.

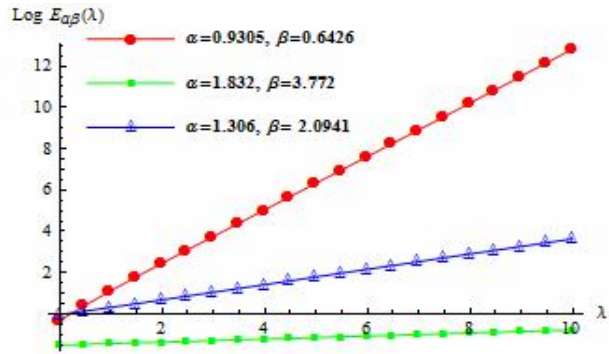

Fig. 11

## 2.7 Sufficient Statistics for $(\alpha, \beta)$ when $\lambda$ is fixed

From equation (5) using the factorization theorem it can easily be seen that for a sample $(x_1, x_2, \cdots, x_n)$ of size $n$, from MLFD$(\lambda, \alpha, \beta)$, $\sum_{i=1}^{n} x_i$ is sufficient statistic for the parameter $\lambda$ when $\alpha, \beta$ are known or fixed.

## 3. MLFD $(\lambda, \alpha, 1)$ as a distribution in a queuing system

MLFD$(\lambda, \alpha, 1)$, like the COM-Poisson distribution [9], can be derived as the probability of the system being in the *k*-th state for a queuing system with state dependent service rate.

Consider a queuing system with Poisson inter arrival times with parameter $\lambda$, first come first serve policy, and exponential service times that depend on the system state (*n*th state means $n$ number of units in the system). The mean service time in the *n*th state is $\mu_n = \mu (n\alpha)_{[\alpha]}, n \geq 1$ and $\mu_n = \mu$ for $\alpha = 0$, where, $1/\mu$ is the normal mean service time for a unit when that unit is the only one in the system and $\alpha$ is the pressure coefficient, a constant reflecting the degree to which the service rate of the system is affected by the system. For the sake of completeness, the proof that the probability is the MLFD$(\lambda, \alpha, 1)$ pmf is given as follows:

Following Conway and Maxwell (p. 134-35 [9]), the system differential difference equations are
$$P_0(t+\Delta) = (1 - \lambda \Delta)P_0(t) + \mu_1 \Delta P_1(t) \tag{8}$$
and $P_n(t+\Delta) = (1 - \lambda \Delta - (n\alpha)_{[\alpha]}\mu \Delta)P_n(t) + \lambda \Delta P_{n-1}(t) + ((n+1)\alpha)_{[\alpha]}\mu \Delta P_{n+1}(t)), n = 1, 2, \cdots$ (9)

From (8) we get $P_0(t+\Delta) - P_0(t) = -\lambda \Delta P_0(t) + \mu \alpha! \Delta P_1(t)$ since $\mu_1 = (\alpha)_\alpha = \alpha!$. This implies,
$$\lim_{\Delta \to 0} \frac{P_0(t+\Delta) - P_0(t)}{\Delta} = -\lambda P_0(t) + \mu \alpha! P_1(t)$$
or $P_0^{/}(t) = -\lambda P_0(t) + \mu \alpha! P_1(t)$.

Assuming a steady state (i.e. $P_n^{/}(t) = 0$ for all *n*) we get
$$P_1(t) = zP_0(t)/\alpha!$$
where $\lambda / \mu = z$. Similarly, from (9), we get
$$\lim_{\Delta \to 0} \frac{P_n(t+\Delta) - P_n(t)}{\Delta} = -(\lambda \Delta - (n\alpha)_{[\alpha]}\mu)P_n(t) + \lambda P_{n-1}(t) + ((n+1)\alpha)_{[\alpha]}\mu P_{n+1}(t))$$

It follows that
$$P_n^{/}(t) = -(\lambda - (n\alpha)_{[\alpha]}\mu)P_n(t) + \lambda P_{n-1}(t) + ((n+1)\alpha)_{[\alpha]}\mu P_{n+1}(t))$$
$$0 = -\mu(z - (n\alpha)_{[\alpha]})P_n(t) + \mu z P_{n-1}(t) + \mu((n+1)\alpha)_{[\alpha]}P_{n+1}(t)),$$
since $P_n^{/}(t) = 0$ for all *n*.
This implies $(z + (n\alpha)_{[\alpha]})P_n(t) = z P_{n-1}(t) + ((n+1)\alpha)_{[\alpha]}P_{n+1}(t))$ since $\mu \neq 0$
Putting $n = 1$ we get
$$(z + \alpha!)P_1(t) = z P_0(t) + (2\alpha)_{[\alpha]}P_2(t)),$$
$$P_2(t) = \{z^2 / (2\alpha)!\}P_0(t). \text{ Since } \alpha!(2\alpha)_{[\alpha]} = (2\alpha)!$$
Similarly, for $n = 2$ we get

$$(z+(2\alpha)_{[\alpha]})P_2(t) = z\ P_1(t) + (3\alpha)_{[\alpha]}P_3(t))$$

$$P_2(t) = \{z^3/(3\alpha)!\}P_0(t),$$

since $(2\alpha!)(3\alpha)_{[\alpha]} = (3\alpha)!$

In general,

$$P_n(t) = \{z^n/(n\alpha)!\}P_0(t),$$

where $P_0(t) = 1/\sum_{n=0}^{\infty}\{z^n/(n\alpha)!\}$. This is the MLFD$(\lambda,\alpha,\beta)$ pmf.

## 4. Reliability and stochastic ordering

Discrete life time model has of late been favourite subject matter of many studies as many times the life of a system is may be measured by counting and even when the life is measured in a continuous scale the actual observations may be recorded in a way making discrete model more appropriate. It is therefore important to study the reliability properties of the proposed discrete distribution. Stochastic ordering is a closely related important area that has found applications in many diverse areas such as economics, reliability, survival analysis, insurance, finance, actuarial and management sciences (see Shaked and Shanthikumar, [58]). In this section we study the reliability properties and stochastic ordering of the MLFD$(\lambda,\alpha,\beta)$ distribution.

***Survival function*:**

$$S(t) = P(X > t) = \lambda^{t+1} E_{\alpha,\beta+(t+1)\alpha}(\lambda)/E_{\alpha,\beta}(\lambda)$$

***Failure rate function***:

$$r(t) = P(X=t)/P(X \geq t) = \{\lambda \Gamma(\alpha t + \beta)E_{\alpha,\beta+(t+1)\alpha}(\lambda)\}^{-1}$$

***Log-concavity***:

The MLFD$(\lambda,\alpha,\beta)$ *has a log-concave probability mass function* since for this distribution (Gupta et al., 1997)

$$\Delta\eta(t) = \frac{P(t+1)}{P(t)} - \frac{P(t+2)}{P(t+1)} = p\frac{\Gamma(\alpha t+\beta)\Gamma(\alpha t+2\alpha\beta)-\{\Gamma(\alpha t+\alpha+\beta)\}^2}{\Gamma(\alpha t+\alpha+\beta)\Gamma(\alpha t+2\alpha\beta)} > 0.$$

The following results are direct consequence of log-concavity (Mark, 1996):

**i.** MLFD$(\lambda,\alpha,\beta)$ *is a **strongly unimodal** distribution* due to the log-concavity of its pmf (see Steutel [64]).

> MLFD$(\lambda,\alpha,\beta)$ ***has a unique mode*** *at* $X = k$ *if*

$$\frac{\Gamma(\alpha k + \beta)}{\Gamma(\alpha k - \alpha + \beta)} < \lambda < \frac{\Gamma(\alpha k + \alpha + \beta)}{\Gamma(\alpha k + \beta)}$$

*Proof*: Since

$P(X = k+1) < P(X = k)$ if $\lambda < \dfrac{\Gamma(\alpha k + \alpha + \beta)}{\Gamma(\alpha k + \beta)}$ and

$P(X = k) > P(X = k-1)$ if $\lambda > \dfrac{\Gamma(\alpha k + \beta)}{\Gamma(\alpha k - \alpha + \beta)}$, the result follows.

- MLFD $(\lambda, \alpha, \beta)$ *has a **non increasing** pmf with a unique mode at* $X = 0$ *if* $\lambda < \frac{\Gamma(\alpha + \beta)}{\Gamma(\beta)}$ (See the pmf plots in Fig. 7 for some choices of $(\lambda, \alpha, \beta)$ satisfying the condition)
- MLFD $(\lambda, \alpha, \beta)$ ***has two modes*** *at* $X = k$ *and* $X = k+1$ *if*
$$\lambda = \Gamma(\alpha k + \alpha + \beta)/\Gamma(\alpha k + \beta)$$
(See the pmf plotted in Fig. 8 for some choices of $(\lambda, \alpha, \beta)$ satisfying the condition)

**ii.** MLFD $(\lambda, \alpha, \beta)$ *is an increasing failure rate function.*
**iii.** MLFD $(\lambda, \alpha, \beta)$ *has at most an exponential tail.*
**iv.** MLFD $(\lambda, \alpha, \beta)$ *remains log concave if truncated.*
**v.** *Convolution of* MLFD $(\lambda, \alpha, \beta)$ *with any other discrete distribution will also result in log concave distribution.*
**vi.** $\frac{P(X = i + k)}{P(X = i)} \geq \frac{P(X = j + k)}{P(X = j)}$ for $i < j$.

*Stochastic ordering with HP*:
The following result stochastically compares the MLFD $(\lambda, \alpha, \beta)$ with the HP $(\lambda, \beta)$ by using the likelihood ratio order.

**Theorem 2.** *For* $\alpha > 1$, $X \sim$ MLFD $(\lambda, \alpha, \beta)$ *is smaller than* HP $(\lambda, \beta)$ *distribution in the likelihood ratio order i.e.* $X \leq_{lr} Y$.

*Proof*: If $X \sim$ MLFD $(\lambda, \alpha, \beta)$ and $Y \sim$ HP $(\lambda, \beta)$ then
$$\frac{P(Y = n)}{P(X = n)} = \frac{\Gamma(n\alpha + \beta)}{\Gamma(n + \beta)} \frac{E_{1,\beta}(\lambda)}{E_{\alpha,\beta}(\lambda)}$$

This is clearly increasing in *n* as $\alpha > 1$ (see Shaked and Shanthikumar [58] and Gupta et al., [19]). Hence the result is proved.

**Remark 3.** For $0 < \alpha < 1$, MLFD $(\lambda, \alpha, \beta)$ *is greater than* HP $(\lambda, \beta)$ *distribution in the likelihood ratio order i.e.* $Y \leq_{lr} X$.

**Corollary:** Again $X \leq_{lr} Y$ implies $X$ is smaller than $Y$ in the hazard rate order and subsequently in the mean residual (MRL) life order (see Gupta et al.[19]). Symbolically,
$$X \leq_{lr} Y \Rightarrow X \leq_{hr} Y \Rightarrow X \leq_{MRL} Y, \text{ for } \alpha > 1.$$

**Corollary:** Since $Y \leq_{lr} X$ implies $Y$ is smaller than $X$ in the hazard rate order and subsequently in the mean residual (MRL) life order (see Gupta et al. [19]). Symbolically,
$$Y \leq_{lr} X \Rightarrow Y \leq_{hr} X \Rightarrow Y \leq_{MRL} X, \text{ for } 0 < \alpha < 1.$$

## 5. Computation of the generalized Mittag-Leffler function $E_{\alpha,\beta}(z)$

For statistical inference and applications it is necessary to compute the generalized Mittag-Leffler function $E_{\alpha,\beta}(z)$, which is the normalizing constant, where

$$E_{\alpha,\beta}(z) = \sum_{k=0}^{\infty} z^k / \{\Gamma(\alpha k + \beta).$$

Numerical computation of the generalized Mittag-Leffler function is well-researched. Seybold and Hilfer [57] gave a numerical algorithm for calculating the generalized Mittag-Leffler function for arbitrary complex argument $z$ and real parameters $\alpha$ and $\beta$ based on a Taylor series, exponentially improved asymptotics and integral representation. If $|z| \leq 1$, a simple way to compute $E_{\alpha,\beta}(z)$ is to calculate the terms $a_k = z^k / \{\Gamma(\alpha k + \beta)$ in the infinite series $E_{\alpha,\beta}(z)$ by employing the recurrence relation

$$a_{k+1} / a_k = z / (\alpha k + \beta)$$

with $a_0 = 1 / \Gamma(\beta)$ (see Lee et al., [38]). This avoids the computation of the gamma function $\Gamma(x)$. The summation is terminated when $a_k$ is very small. The error estimate is given by Theorem 4.1 of Seybold and Hilfer [57] which determines the number of terms $N$ such that

$$E_{\alpha,\beta}(z) \approx \sum_{k=0}^{N} z^k / \{\Gamma(\alpha k + \beta).$$

For other values of z, asymptotic series and integral representation (Seybold and Hilfer [57], equations (2.3), (2.4) and (2.7) respectively) are employed. Error estimates are also given for these cases. The computation of the Mittag-Leffler function is given by many software packages like Matlab (MLF (alpha, Z, P)) and Mathematica ( MittagLefflerE[*a*, *b*, *z*]). (See also Gorenflo et al., [17]).

## 6. Data analysis with MLFD $(\lambda, \alpha, \beta)$

### 6.1 Parameter estimation

Suppose that we have a sample of size *n* from MLFD $(\lambda, \alpha, \beta)$ reported in grouped frequency in *k* classes, like $(\mathbf{X}, \mathbf{f}) = \{(x_1, f_1), (x_2, f_2), \ldots, (x_k, f_k)\}$, where $f_i$ is the frequency of the $i^{th}$ observed value $x_i$ and $n = \sum_{i=1}^{k} f_i$ is the sample size. Then the log-likelihood function is given by

$$\ln L(x_1, x_2, \cdots, x_k | \lambda, \alpha, \beta) = (\sum_{i=1}^{k} f_i x_i) \ln p - \sum_{i=1}^{k} f_i \ln \Gamma(\alpha x_i + \beta) - n \ln \sum_{j=0}^{\infty} \{p^j / \Gamma(\alpha x_j + \beta)\}.$$

Numerical optimization method is used to obtain the maximum likelihood estimates of the parameters required for the data fitting experiments and the likelihood ratio test.

The estimates of the variances and covariances of the mles can be obtained from the inverse of the Fishers information matrix.

### 6.2 Likelihood ratio test

The MLFD $(\lambda, \alpha, \beta)$ reduces to the HP distribution with parameters $(\lambda, \beta)$ when $\alpha = 1$ and MLFD $(\lambda, \alpha, 1)$ when $\beta = 1$. Since the distributions are nested we have employed the likelihood ratio criterion to test the following null hypotheses:

**I.** HP $(\lambda, \beta)$ vs. MLFD $(\lambda, \alpha, \beta)$:

$H_0: \alpha = 1$, that is the sample is from HP $(\lambda, \beta)$ distribution against the alternative

$H_1: \alpha \neq 1$, that is the sample is from MLFD$(\lambda, \alpha, \beta)$.

**II.** MLFD$(\lambda, \alpha, 1)$ vs. MLFD$(\lambda, \alpha, \beta)$: Here the null hypothesis is:

$H_0: \beta = 1$, that is the sample is from MLFD$(\lambda, \alpha, 1)$ against the alternative

$H_1: \beta \neq 1$, that is the sample is from MLFD$(\lambda, \alpha, \beta)$.

Writing $\boldsymbol{\theta} = (\lambda, \alpha, \beta)$ the likelihood ratio test statistic is given by LR $= -2\ln L(\hat{\boldsymbol{\theta}}^*; x)/L(\hat{\boldsymbol{\theta}}, x)$, where $\hat{\boldsymbol{\theta}}^*$ is the restricted ML estimates under the null hypothesis $H_0$ and $\hat{\boldsymbol{\theta}}$ is the unrestricted ML estimates under the alternative hypothesis $H_1$. Under the null hypothesis $H_0$ the LR criterion follows Chi-square distribution with 1 degree of freedom.

*6.3 Numerical examples*

Three data sets from the literature have been considered here. The first set is the frequency data (Table 1.1) on insurance claims and incapacity caused by sickness or accident (Lundberg [40] also used by Phang, et. al. [48]) is over dispersed with index of dispersion 2.28255, the second set (Table 2.1) is the frequency data on the number of sickness absences (1955–1964) Taylor [65] is an over dispersed data with index of dispersion 5.38588, and the third set (Table 3.1) is the frequency data on the secondary association of chromosomes in Brassika (Skellam, [61]) is an under dispersed date with index of dispersion 0.49155.

The MLFD$(\lambda, \alpha, \beta)$ model has been fitted and compared with two related distributions namely the HP$(\lambda, \beta)$, MLFD$(\lambda, \alpha, 1)$ distributions. Asymptotic variance and covariance of the MLEs are estimated from the inverse of the information matrix. The performances of various distributions are compared using the log likely hood value and the AIC (Akaike Information Criterion) defined as AIC $= -2 \ln L + 2k$, where $k$ is the number of parameter(s) and $\ln L$ is the maximum of log-likelihood for a given data set (Burnham and Anderson [4]). Relatively high (small) values of log likelihood (AIC) would indicate better fitting. As seen in Tables 1.2, 2.2 and 2.3, MLFD$(\lambda, \alpha, \beta)$ clearly gives the best fit in all the three cases.

In addition the likelihood ratio test is also performed to discriminate between the MLFD$(\lambda, \alpha, \beta)$ and its sub models namely HP$(\lambda, \beta)$ and MLFD$(\lambda, \alpha, 1)$. It is seen from Tables 1.3, 2.3 and 3.3 that for the three examples MLFD$(\lambda, \alpha, \beta)$ is chosen over HP$(\lambda, \beta)$. MLFD$(\lambda, \alpha, \beta)$ is also chosen over MLFD$(\lambda, \alpha, 1)$ in example 2 and 3.

**Table 1.1**
Data on insurance claims and incapacity caused by sickness or accident [40]

| X | 0 | 1 | 2 | 3 | 4 | 5 | 6 | 7 | 8 | 9 | 10 | 11 | 12 | 13 | 14 | 15 | Total |
|---|---|---|---|---|---|---|---|---|---|---|----|----|----|----|----|----|-------|
| Freq | 187 | 185 | 200 | 164 | 107 | 68 | 49 | 39 | 21 | 12 | 11 | 2 | 5 | 2 | 3 | 1 | 1056 |

Mean = 2.80587, Variance = 6.40455, ID = **2.28255**.

**Table 1.2**
Summary of results of fittings of the data set in Table 1.1

| Distribution | MLEs | Log L | AIC |
|---|---|---|---|
| HP$(\lambda, \beta)$ | $(\lambda, \beta) = (13.0675, 13.6849)$ | -2270.16 | 4544.32 |
| MLFD$(\lambda, \alpha, 1)$ | $(\lambda, \alpha) = (0.9128, 0.3332)$ | -2267.01 | 4538.03 |
| MLFD$(\lambda, \alpha, \beta)$ | $(\lambda, \alpha, \beta) = (0.6114, 0.1282, 0.1603)$ | **-2265.39** | **4536.78** |

**Table 1.3**
Summary of results LR tests of data set in Table 1.1

| $H_0$ | $H_1$ | LR Statistic | p-value |
|---|---|---|---|
| $\alpha = 1$ (HP) | $\alpha \neq 1$ (MLFD $(\lambda, \alpha, \beta)$) | 9.5451 | 0.0020 |
| $\beta = 1$ (MLFD $(\lambda, \alpha, 1)$) | $\beta \neq 1$ (MLFD $(\lambda, \alpha, \beta)$) | 3.2490 | 0.0715 |

Estimated variances and covariance's of the mles of MLFD $(\lambda, \alpha, \beta)$:
$Var(\hat{\lambda}) = 0.0036$, $Var(\hat{\alpha}) = 0.0078$, $Var(\hat{\beta}) = 0.0257$, $Cov(\hat{\lambda}, \hat{\alpha}) = 0.0052$
$Cov(\hat{\lambda}, \hat{\beta}) = 0.0096$, $Cov(\hat{\alpha}, \hat{\beta}) = 0.01394$.

**Table 2.1**
Taylor's [65] data on the number of sickness absences (1955–1964)

| x | 0 | 1 | 2 | 3 | 4 | 5 | 6 | 7 | 8 | 9 | 10 | 11 | 12 | 13 | |
|---|---|---|---|---|---|---|---|---|---|---|---|---|---|---|---|
| Freq | 31 | 35 | 55 | 59 | 49 | 59 | 41 | 38 | 32 | 31 | 24 | 22 | 22 | 17 | |
| x | 14 | 15 | 16 | 17 | 18 | 19 | 20 | 21 | 22 | 23 | 24 | 25 | 26 | 27 | Total |
| Freq | 16 | 10 | 8 | 8 | 11 | 4 | 3 | 6 | 9 | 4 | 7 | 6 | 8 | 8 | 623 |

Mean = 7.95666, Variance = 42.8537, ID = **5.38588**

**Table 2.2**
Summary of results of fittings of the data set in Table 2.1

| Distribution | MLEs | Log L | AIC |
|---|---|---|---|
| HP | $(\lambda, \beta) = (113.896, 115.891)$ | -1927.92 | 3859.83 |
| MLFD $(\lambda, \alpha, 1)$ | $(\lambda, \alpha) = (0.975, 0.1319)$ | -1925.52 | 3855.04 |
| MLFD $(\lambda, \alpha, \beta)$ | $(\lambda, \alpha, \beta) = (0.8057, 0.2989 \times 10^{-4}, 0.4338 \times 10^{-4})$ | **-1920.42** | **3846.84** |

Estimated variances and covariance's of the mles of MLFD $(\lambda, \alpha, \beta)$:
$Var(\hat{\lambda}) = 0.0006$, $Var(\hat{\alpha}) = 0.0019$, $Var(\hat{\beta}) = 0.0040$, $Cov(\hat{\lambda}, \hat{\alpha}) = -0.0010$,
$Cov(\hat{\lambda}, \hat{\beta}) = -0.0014$, $Cov(\hat{\alpha}, \hat{\beta}) = 0.0027$.

**Table 2.3**
Summary of results LR tests of data set in Table 2.1

| $H_0$ | $H_1$ | LR Statistic | p-value |
|---|---|---|---|
| $\alpha = 1$ (HP) | $\alpha \neq 1$ (MLFD $(\lambda, \alpha, \beta)$) | 14.9919 | 0.000108 |
| $\beta = 1$ (MLFD $(\lambda, \alpha, 1)$) | $\beta \neq 1$ (MLFD $(\lambda, \alpha, \beta)$) | 10.1998 | 0.001405 |

**Table 3.1**
Dataset consists of the secondary association of chromosomes in Brassika (Skellam, [61])

| X | 0 | 1 | 2 | 3 | Total |
|---|---|---|---|---|---|
| Freq | 32 | 103 | 122 | 80 | 337 |

Mean = 1.74184, Variance = 0.856202, ID = 0.49155.

**Table 3.2** Summary of results of fittings of the data set in Table 3.1

| Distribution | MLEs | Log L | AIC |
|---|---|---|---|
| HP | $(\lambda, \beta) = (1.0409, 0.2337)$ | -468.636 | 941.272 |
| $MLFD(\lambda, \alpha, 1)$ | $(\lambda, \alpha) = (8.9113, 1.8574)$ | -460.317 | 924.635 |
| $MLFD(\lambda, \alpha, \beta)$ | $(\lambda, \alpha, \beta) = (36.6589, 2.2862, 1.9074)$ | **-458.208** | **922.417** |

Estimated variances and covariance's of the mles of $MLFD(\lambda, \alpha, \beta)$:
$Var(\hat{z}) = 139.428$, $Var(\hat{\alpha}) = 0.0106$, $Var(\hat{\beta}) = 0.1083$, $Cov(\hat{z}, \hat{\alpha}) = 1.0694$, $Cov(\hat{z}, \hat{\beta}) = 2.6506$, $Cov(\hat{\alpha}, \hat{\beta}) = 0.0098$.

**Table 3.3** Summary of results LR tests of data set in Table 3.1

| $H_0$ | $H_1$ | LR Statistic | p-value |
|---|---|---|---|
| $\alpha = 1$ (HP) | $\alpha \neq 1$ ($MLFD(\lambda, \alpha, \beta)$) | 20.855 | 0.005 |
| $\beta = 1$ ($MLFD(\lambda, \alpha, 1)$) | $\beta \neq 1$ ($MLFD(\lambda, \alpha, \beta)$) | 4.21799 | 0.040 |

## 7. Conclusion

A new generalization of the hyper-Poisson distribution which is a continuous bridge between geometric and hyper-Poisson is derived using the generalized Mittag Leffler function. Some known and new distributions are seen as particular cases of this distribution. This new generalization belongs to the generalized power series, generalized hyper geometric families and also arises as weighted Poisson distributions. Like the hyper-Poisson, COM-Poisson and generalized Poisson distributions, this distribution is also able to model under, equi- and over dispersion. Although the new generalization of the hyper-Poisson distribution has an extra parameter it is computational not more complicated than the hyper-Poisson since it retains the two-term probability recurrence formula and the normalizing constant, in terms of the generalized Mittag Leffler function is readily computed. It has many interesting probabilistic and reliability properties and is found to be a better empirical model than the hyper-Poisson distribution.

*Acknowledgement*: The authors wish to acknowledge support in parts from the Ministry of Education FRGS grant FP010-2013A and University of Malaya's UMRGS grant RP009A-13AFR.

**References**
1. M. Ahmad, A short note on Conway-Maxwell-Hyper Poisson Distribution, Pak. J. Stat. 23 (2007), 135–137.
2. G. Antic, E. Stadlober, P. Grzybek, E. Kelih, Word Length and Frequency Distributions in Different Text Genres. In: Spiliopoulou, Myra; Kruse, Rudolf; Nürnberger, Andreas; Borgelt, Christian; Gaul, Wolfgang (eds.), From Data and Information Analysisto Knowledge Engineering. Heidelberg, Berlin: Springer, (2006) 310–317.

3. K. Best, Kommentierte Bibliographie zum Göttinger Projekt. In: Best, Karl-Heinz (ed.), Häufigkeitsverteilungen in Texten. Göttingen: Pest & Gutschmidt, (2001) 248–310.

4. K.P. Burnham, D.R. Anderson, Multimodel Inference, Understanding AIC and BIC in model selection, Sociological Methods and Research, 33 (2004) 261304

5. G.E. Bardwell, E.L. Crow, A two Parameter Family of Hyper- Poisson Distributions, Journal of the American Statistical Association, 59 (1964) 133-141, doi: 10.1080/01621459.1964.10480706.

6. A.C. Cameron, P. Johansson, Count Data Regression Using Series Expansions: with Applications, Journal of Applied Econometrics, 12, 3 (1997) 203–223.

7. J.D. Castillo, M. Pérez-Casany, Over-dispersed and under-dispersed Poisson Generalizations, Journal of Statistical Planning and Inference, 134 (2005) 486-500.

8. P.C. Consul, Generalized Poisson Distributions: Properties and Applications, Marcel Dekker Inc, New York / Basel, (1989).

9. R.W. Conway, W.L. Maxwell, A queuing model with state dependent service rates, Journal of Industrial Engineering, 12 (1962) 132–136

10. E.L. Crow, G.E. Bardwell Estimation of the parameters of the Hyper-Poisson Distributions. In: Patil, G.P. (ed) Classical and Contagious Discrete Distributions, 27–140, Pergamon, Oxford, (1965).

11. L. Devroye, A note on Linnik's Distribution, Statistics and Probability Letters, 9 (1990) 305–306.

12. B. Efron, Double Exponential-Families and their use in Generalized Linear-Regression. Journal of the American Statistical Association, 81, 395(1986) 709–721.

13. A. Erdelyi, Higher Transcendental Functions, Vol. 3, McGrwa-Hill, New York, (1955).

14. B. Fisher, A. Kilicman, Some Results on the Gamma Function for Negative Integers, Appl. Math. Inf. Sci. 6, 2(2012) 173-176.

15. R.A. Francis, S.R. Geedipally, S.D. Guikema, S.S Dhavala, D. Lord, S. LaRocca, Characterizing the Performance of the Conway–Maxwell Poisson Generalized Linear Model. Risk Analysis, 32, 1, (2012) 167–183.

16. S. Gerhold, Asymptotics for a variant of the Mittag-Lefller Function, Integral Transforms and Special functions, 23(6) (2012) 397–403.

17. R. Gorenflo, L. Joulia, Y. Luchko, (2002) Computation of the Mittag Leffler Function $E_{\alpha,\beta}(z)$ and its derivatives, Fractional Calculus & Applied Analysis, 5 (4) 491-518.

18. P.L. Gupta, R.C. Gupta, R.C. Tripathi, On the Monotonic properties of Discrete Failure Rates. Journal of Statistical Planning and Inference, 65(1997) 255–268

19. R.C. Gupta, S.Z. Sim, S.H. Ong, Analysis of discrete data by Conway-Maxwell Poisson distribution. AStA Advances in Data Analysis (2014), http://dx.doi.org/10.1007/s10182-014-0226-4.

20. J.W. Hanneken, B.N.N. Achar, R. Puzio, D.M. Vaught, Properties of the Mittag–Leffler function for negative alpha, Phys. Scr. T136 (2009), 014037 (doi:10.1088/0031-8949/2009/T136/014037), 1-5.

21. H.J. Haubold, A.M. Mathai, R.K. Saxena, Mittag-Leffler Functions and their Applications, Journal of Applied Mathematics (2011), Article ID 298628, 51 pages,doi:10.1155/2011/298628

22. K.G. Janardan, (1984) Moments of certain series distributions and their applications, Siam Journal of Applied Mathematics, 44, 4, 854-868.

23. K. Jayakumar, R.N. Pillai, The first order autoregressive Mittag-Leffler process, Journal of Applied Probability, 30 (1993), 462-466.

24. K. Jayakumar, On Mittag-Leffler process, Mathematical and Computer Modelling, 37 (2003) 1427-1434.

25. N.L. Johnson, A.W. Kemp, S. Kotz, Univariate Discrete Distributions. Wiley, New York (2005).

26. K.K. Jose, R.N. Pillai, Generalized autoregressive time series models in Mittag-Leffler variables, Recent Advances in Statistics (1996) 96-103.

27. K.K. Jose, P. Uma, V. Seethalekshmi, H.J. Haubold, (2010) Generalized Mittag-Leffler processes for applications in astrophysics and time series modelling, Astrophysics and Space Science Proceedings, doi: 10.1007/978-3-642-03325-4.

28. K.K. Jose, B. Abraham, A count data model based on Mittag-Leffler inter arrival times, Statistica, anno LXXI, 4 (2011) 501-514.

29. C.D. Kemp, q-analogues of the Hyper-Poisson Distribution. Journal of Statistical Planning and Inference, 101 (2002) 179–183

30. C.D. Kemp, A.W. Kemp, Generalized hypergeometric distributions. Journal of the Royal Statistical Society, Series B (Methodological) 18(2) (1956) 202-211.

31. S.H. Khazraee, A.J. Sáez-Castillo, S.R. Geedipally, D. Lord, Application of the Hyper-Poisson generalized linear model for analyzing motor vehicle crashes- 93rd Annual Meeting of the Transportation Research Board (2013).

32. S.K. Kattumannit, L. Tibiletti, Moment identity for discrete random Variable and its applications, Statistics: A Journal of Theoretical and Applied Statistics, doi: 10.1080/02331888.2011.555548 (2011).

33. C.G. Khatri, On certain properties of power-series distributions, Biometrika, 46 (1959), 486-490.

34. C.S. Kumar, B.U. Nair, A Modified version of Hyper-Poisson Distribution and its applications. J. Stat. Appl. 6 (2011) 23-34

35. C.S. Kumar, B.U. Nair, An alternative Hyper-Poisson Distribution. Statistica 3 (2012), 357-369

36. C.S. Kumar, B.U. Nair, Modified alternative Hyper-Poisson Distribution. In: Kumar, C.S., Chacko, M., Sathar, E.I.A. (eds) Collection of Recent Statistical Methods and Applications, Department of Statistics, University of Kerala publication, Trivandrum (2013) 97-109..

37. C.S. Kumar, B.U. Nair, On a class of Hyper-Poisson and alternative Hyper-Poisson Distributions, OPSEARCHI (Operational Research Society of India) (2014), doi: 10.1007/s12597-013-0169-7

38. P.A. Lee, S.H. Ong, H.M. Srivastava, Some integrals of the products of Laguerre Polynomials. Internat J Comput Math, 78 (2001) 303–321

39. G.D. Lin, On Mittag-Leffler Distribution, Journal of Statistical Planning and Inference, 74 (1998) 1 - 9.

40. O. Lundberg, On Random Processes and Their Application to Sickness and Accident Statistics. Uppsala: Almquist and Wiksells (1940).

41. Y.A. Mark, Log-concave Probability Distributions: Theory and Statistical testing. Working paper, 96-01 (1996), Published by Center for Labour Market and Social Research, University of Aarhus and the Aarhus School of Business.

42. T.P. Minka, , G. Shmueli, J.B. Kadane, S. Borle, P. Boatwright, Computing with the COM-Poisson Distribution. Technical Report #776 (2003), Department of Statistics, Carnegie Mellon University. http://www.stat.cmu.edu/tr/tr776/tr776.html.

43. G. Mittag-Leffler, Sur la nouvelle fonction $E_\alpha(x)$, Comptes Rendus Acad. Sci. Paris 137 (1903) 554-558.

44. G. Mittag-Leffler, Sur la representation analytique d'une branche uniforme d'une fonction monogene, Acta Math. 29 (1905) 101-181.

45. S.H. Ong, S. Chakraborty, T. Imoto, K. Shimizu Generalization of non-central negative binomial, Charlier series distributions and their extensions by Lagrange expansion. Communications in Statistics-Theory and Methods 41(4) (2012) 571–590.

46. G.P. Patil, Estimation by two moments method for generalized power series distribution and certain applications, Sankhya, B, 24, 3&4 (1962) 201-214.

47. G.P. Patil, Estimation for the generalized power series distribution with two parameters and its application to binomial distribution, Contributions to Statistics, C.R. Rao (editor), 335–344. Calcutta: Statistical Publishing Society; Oxford: Pergamon, (1964).

48. Y.N. Phang, S.Z.Sim, S.H. Ong, Statistical Analysis for the Inverse Trinomial Distribution, Communications in Statistics - Simulation and Computation, 42, 9 (2013) 2073-2085.

49. G.P. Patil, C.R. Rao, M.V. Ratnaparkhi, On discrete weighted distributions and their use in model for observed data. Communications in Statistics-Theory and Methods, 15 (3) (1986) 907-918.

50. R.N. Pillai, On Mittag Lefller functions and related distributions. Ann. Inst. Statist. Math. 42 (1990) 157-161.

51. R.N. Pillai, K. Jayakumar, Discrete Mittag-Leffler Distributions, Statistics & Probability Letters, 23 (1995) 271-274.

52. T.R. Prabhakar, A singular integral equation with a generalized Mittag-Leffler function in the Kernel, Yokohama Math. J. No 19 (1971) 7-15.

53. A. Roohi, M. Ahmad, Estimation of the parameter of hyper-Poisson distribution using negative moments. Pak. J. Stat. 19 (2003a) 99–105.

54. A. Roohi, M. Ahmad, Inverse ascending factorial moments of the hyper-Poisson probability distribution. Pak. J. Stat. 19 (2003b) 273–280

55. S.M. Ross, Stochastic Processes, Wiley, NewYork (1983).

56. A.J. Sáez-Castillo, A. Conde-Sánchez, A Hyper-Poisson regression model for Overdispersed and Underdispersed Count Data. Computational Statistics and Data Analysis (2012), doi:10.1016/j.csda.2012.12.009.

57. H. Seybold, R. Hilfer, Numerical algorithm for calculating the generalized Mittag - Leffler function. SIAM Journal on Numerical Analysis, 47(1) (2008) 69-88.

58. M. Shaked, J. G. Shanthikumar, Stochastic Orders. Springer Verlag (2007).

59. G. Shmueli, T.P. Minka, J.B. Kadane, S. Borle, S. Boatwright, A useful distribution for fitting discrete data-revival of the Conway-Maxwell Poisson distribution. Applied Statistics, 54 (2005)127–142

60. T. Simon, Mittag - Leffler functions and complete monotonicity, arXiv: 1312. 4513v [math.CA] (2013).

61. J.G. Skellam, A probability distribution derived from the Binomial distribution by regarding the probability of success as variable between the sets of trials, Journal of the Royal Statistical Society, Series B, 10 (1948) 257-261.

62. P.J. Staff, The displaced Poisson distribution, Australian Journal of Statistics, 6 (1964) 12-20.

63. P.J. Staff, The displaced Poisson distribution. Region B, Journal of American Statistical Association, 62, 318 (1967) 643-654.

64. F.W. Steutel, Log-concave and log-convex distributions. Volume 5, Encyclopedia of Statistical Sciences (eds S. Kotz, N. L. Johnson and C. B. Read) (1985).

65. P.J. Taylor, Individual variations in sickness absence. Brit. J. Industr. Med. 24 (1967) 169.

66. G. Walther, Inference and Modeling with Log-concave Distributions. Statistical Science 24(3) (2009) 319–327.